\documentclass[twoside]{article}
\usepackage{amsmath,bm}
\usepackage{amssymb,amsfonts}
\usepackage{graphicx}                 
\usepackage{float}

\newcommand{\bE}{\mathbf{E}}

\newcommand{\bH}{\mathbf{H}}

\newcommand{\bx}{\mathbf{x}}

\newcommand{\cD}{\mathcal{D}}

\newcommand{\cP}{\mathcal{P}}

\newcommand{\cE}{\mathcal{E}}

\begin{document}


\title{Isoptic surfaces of polyhedra}

\author{G\'eza~Csima and Jen\H o~Szirmai, \\
Budapest University of Technology and Economics, \\
Institute of Mathematics,
Department of Geometry Budapest, \\
P.O. Box 91, H-1521 \\
csgeza@math.bme.hu, szirmai@math.bme.hu}
\maketitle
\thanks{}



\begin{abstract}
The theory of the isoptic curves is widely studied in the Euclidean plane $\bE^2$ (see \cite{CMM91}
and \cite{Wi} and the references given there). The analogous question was investigated by the authors
in the hyperbolic $\bH^2$ and elliptic $\cE^2$ planes (see \cite{CsSz1}, \cite{CsSz2}, \cite{CsSz5}), but in the higher dimensional
spaces there are only a few result in this topic.

In \cite{CsSz4} we gave a natural extension of the notion of the isoptic curves to the
$n$-dimensional Euclidean space $\bE^n$ $(n\ge 3)$
which are called isoptic hypersurfaces. Now we develope an algorithm to determine
the isoptic surface $\mathcal{H}_{\cP}$ of a $3$-dimensional polytop $\mathcal{P}$.

We will determine the isoptic surfaces for Platonic solids and for some semi-regular Archimedean polytopes
and visualize them with Wolfram Mathematica.
\end{abstract}

\newtheorem{theorem}{Theorem}[section]
\newtheorem{corollary}[theorem]{Corollary}
\newtheorem{lemma}[theorem]{Lemma}
\newtheorem{exmple}[theorem]{Example}
\newtheorem{defn}[theorem]{Definition}
\newtheorem{rmrk}[theorem]{Remark}
\newenvironment{definition}{\begin{defn}\normalfont}{\end{defn}}
\newenvironment{remark}{\begin{rmrk}\normalfont}{\end{rmrk}}
\newenvironment{example}{\begin{exmple}\normalfont}{\end{exmple}}
\newtheorem{remarque}{Remark}



\section{Introduction}
Let $G$ be one of the constant curvature plane geometries, either the Euclidean $\bE^2$ or the hyperbolic $\bH^2$ or the elliptic $\cE^2$.
The isoptic curve of a given plane curve $\mathcal{C}$ is the locus of points $P \in G$ where $\mathcal{C}$ is seen under
a given fixed angle $\alpha$ $(0<\alpha <\pi)$.
An isoptic curve formed by the locus of tangents meeting at right angles is called orthoptic curve.
The name isoptic curve was suggested by Taylor in \cite{T}.

In \cite{CMM91} and \cite{CMM96} the isoptic curves of the closed, strictly convex curves are studied, using their support function.
The papers \cite{Wu71-1} and \cite{Wu71-2} deal with curves having a
circle or an ellipse for an isoptic curve. Further curves appearing as isoptic curves are well studied in the Euclidean plane geometry
$\bE^2$, see \textit{e.g.} \cite{Lo}, \cite{Wi}.
Isoptic curves of conic sections have been studied in \cite{H} and \cite{S}.
A lot of papers concentrate on the properties of the isoptics \textit{e.g.} \cite{MM}, \cite{M}, and the reference given there.
The papers \cite{Kur1} and \cite{Kur2} deal with inverse problems.

In the hyperbolic and elliptic planes $\bH^2$ and $\cE^2$ the isoptic curves of segments and proper conic sections are
determined by the authors (\cite{CsSz1}, \cite{CsSz2}, \cite{CsSz3}). In \cite{CsSz5} we extended the notion of the isoptic
curves to the outer (non-proper) points of the hyperbolic plane and determined the isoptic curves of the generalized conic sections.

It is known, that the angle of two half-lines with vertex $A$ in the plane can be measured by the arc length on the unit
circle around the point $A$.
This statement can be generalized to the higher dimensional Euclidean spaces.
The notion of the {\it solid angle} is well known and widely studied in the literatur (see \cite{Gard}).
We recall this definition concerning the $3$-dimensional Euclidean space $\mathbf{E}^3$ and we will apply it in the following.
\begin{defn}
The solid angle $\Omega_S(\mathbf{p})$ subtended by a surface
$S$ is defined as the surface area of a unit sphere covered by the surface's projection onto the sphere around $P(\mathbf{p})$.
\end{defn}
Solid angle is measured in {\it steradians}, e.g. the solid angle corresponding to all of space being subtended is $4\pi$ steradians.
Moreover, this notion has several important applications in physics (particular in astrophysics, radiometry or photometry)
(see \cite{Camp}) even in computational geometry (see \cite{Joe}) and we will use it for the definition of the isoptic surfaces.

The isoptic hypersurface in the $n$-dimensional Euclidean space $(n\ge3)$ is defined in \cite{CsSz4} and
now, we recall some statements and specify them to the $\bE^3$.
\begin{defn}
The isoptic hypersurface $\mathcal{H}_{\cD}^{\alpha}$ in $\mathbf{E}^3$ of an arbitrary $3$-dimensional compact domain $\mathcal{D}$
is the locus of points $P$ where the measure of the projection of $\mathcal{D}$ onto the unit sphere around $P$ is a given fixed value
$\alpha$ $(0<\alpha < 2\pi)$ (see Fig.~1).
\end{defn}

\begin{figure}[ht]
\begin{center}
\includegraphics[width=8cm]{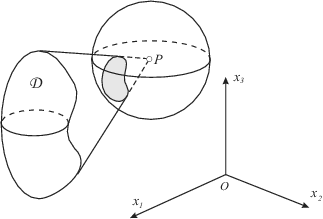}
\caption{Projection of a compact domain $\cD$ to unit sphere in $\bE^3$}
\end{center}
\end{figure}

In this paper we develope an algorithm and the corresponding computer program to determine the isoptic surface of an arbitrary convex
polyhedron in the $3$-dimensional Euclidean space. We apply this algorithm for the regular Platonic solids and some semi-regular Archimedean solids as well. 
This generalization of the isoptic curves to the $3$-dimensional space provide several new way to extend the notion of isoptic surfaces to the "smooth surfaces" and for applications.
\section{The algorithm}
In this section we discuss the algorithm developed to determine the isoptic surface of a given polyhedron.

\begin{enumerate}
\item We assume that an arbitrary polyhedron $\mathcal{P}$ is given by the usual data structure. This consists of the list of facets $\mathcal{F}_\mathcal{P}$ with the set of vertices $V_i$ in clockwise order. Each facet can be embedded into a plane.
%
%

It is well known, that if $\mathbf{a}\in \mathbf{R}^3\setminus \{\mathbf{0}\}$ and $b\in \mathbf{R}$ then
$\{\mathbf{x}\in\mathbf{R}^3|\mathbf{a}^T\mathbf{x} = b \}$ is a plane and $\{\mathbf{x}\in\mathbf{R}^3|\mathbf{a}^T\mathbf{x} \le b \}$ define a halfspace. Every polyhedron is the
intersection of finitely many halfspaces. Therefore an arbitrary polyhedron can also be given by a system of inequality
$A \mathbf{x}\le \mathbf{b}$ where $A\in\mathbf{R}^{m\times3}$ ($4 \leq m\in\mathbb{N}$), $\mathbf{x}\in \mathbf{R}^3$
and $\mathbf{b}\in \mathbf{R}^m$. The solution set of the above equation system provide the polyhedron.
\item 

We have to decide for an arbitrary point $P(\mathbf{p})\in \mathbf{E}^3$, that which faces of $\mathcal{P}$ "can be seen" from it.
Let us denote the $i^{\mathrm{th}}$ facet of $\mathcal{P}$ with $\mathcal{F}_\mathcal{P}^i$ $(i=1\dots m)$ and $\mathbf{a}^i$ is the vector
derived by $i^{\mathrm{th}}$ row of the matrix $A$ which characterize the face $\mathcal{F}_\mathcal{P}^i$.

Since the polyhedron $\mathcal{P}$ is given by the inequality system ${A}\mathbf{x} \leq \mathbf{b}$, where all inequality
${\mathbf{a}^i} \mathbf{x} \le b_i $ ($i\in \{1,2\dots,m\}$) is assigned to a certain face, therefore the facet $\mathcal{F}_\mathcal{P}^i$ is visible from $P$ iff the inequality $\mathbf{a}^i \mathbf{p} > {b_i}$ holds.

Now, we define the characteristic function $\mathbb{I}^{i}_{\mathcal{P}}(\bx)$ for each face $\mathcal{F}_\mathcal{P}^i$:
\[
\mathbb{I}^{i}_{\mathcal{P}}(\bx)=\left\{
\begin{array}{lr}
	1 & \mathbf{a}^i \mathbf{x}> {b_i}       \\
	0 & \mathbf{a}^i \mathbf{x} \leq {b_i}
\end{array}
\right.
\]
\item Using the Definition 1.1, let $\Omega_i(\mathbf{p}):=\Omega_{\mathcal{F}^{i}_{\mathcal{P}}}(\mathbf{p})$ 
be the solid angle of the facet $\mathcal{F}_\mathcal{P}^i$  regarding the point $P(\mathbf{p})$.

To determine $\Omega_i(\mathbf{p})$, we use the machinery of the spherical geometry.
Let us suppose that $\mathcal{F}_\mathcal{P}^i$ contains $n_i$ vertices, $V_{i_j}(\bx_{i_j})$, $(j=1\dots n_i)$
where the vertices are given in clockwise order.
Projecting these vertexes onto the unit sphere around $P(\mathbf{p})$ we get a spherical $n_i$-gon (see Fig~2)
\begin{figure}[ht]
\begin{center}
\includegraphics[width=8cm]{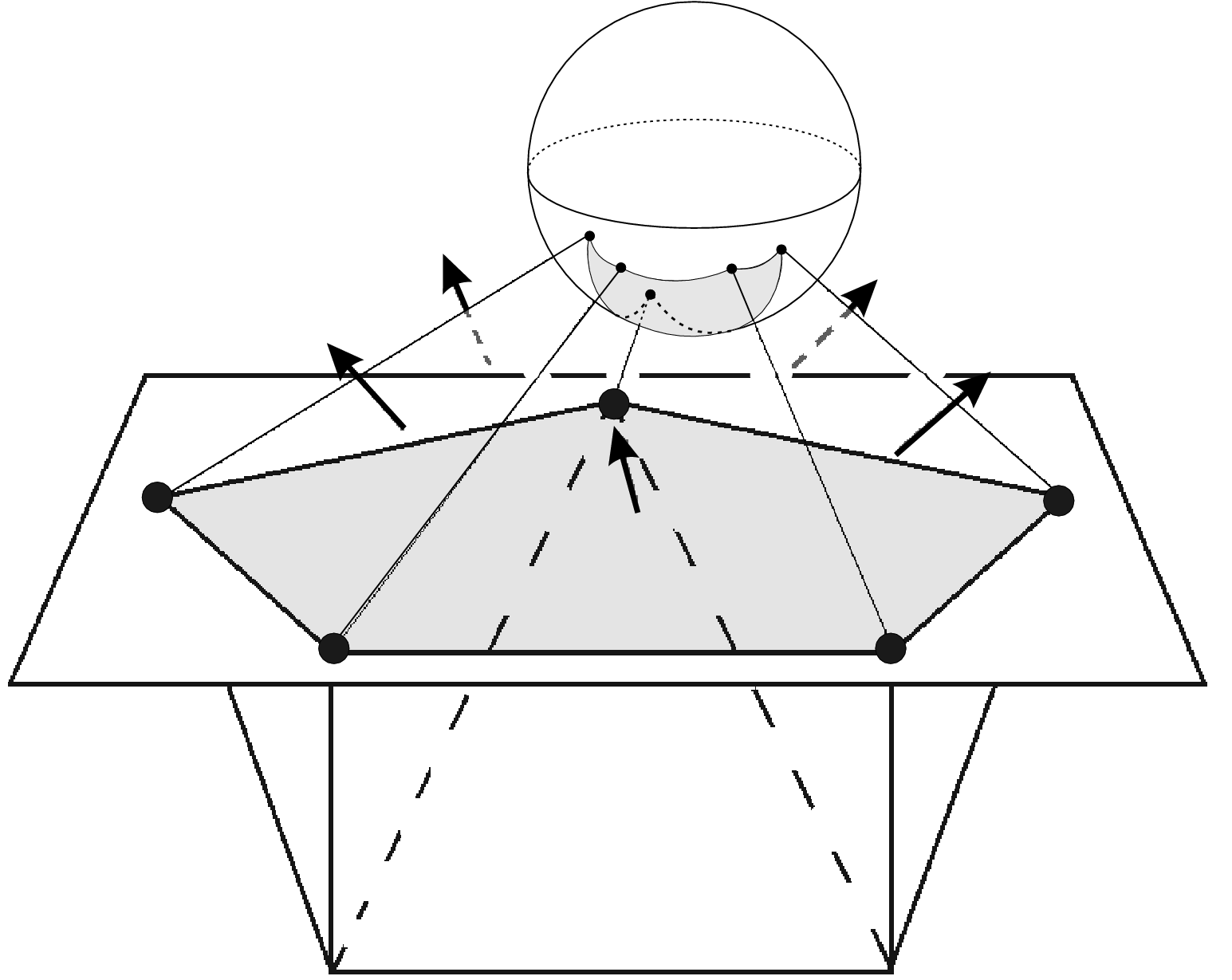}
\caption{Projection of a compact domain $\cD$ to unit sphere in $\bE^3$}
\end{center}
\end{figure}
whose area can be calculated by the usual formula
\[
\Omega_i(\mathbf{p})=(\Theta-(n_i-2)\pi).
\]
Here $\Theta$ is the sum of the angles $\tau_j$ of the spherical polygon $\mathcal{F}_\mathcal{P}^i$ 
where the angles are measured in radians.
\item To obtain these $\tau_j$ angles, we need to determine the angles of the planes that contain
the neigbouring edges $\overline{PV}_{i_k}$ where $k=j-1, j$ and $k=j, j+1$:
\[
\tau_j=\pi-\arccos\left(\frac{\left\langle \overrightarrow{PV}_{i_{j-1}}\times\overrightarrow{PV}_{i_j},\overrightarrow{PV}_{i_{j}}\times\overrightarrow{PV}_{i_{j+1}} 
\right\rangle}{\left|\overrightarrow{PV}_{i_{j-1}}\times\overrightarrow{PV}_{i_j}\right|\left|\overrightarrow{PV}_{i_{j}}
\times\overrightarrow{PV}_{i_{j+1}}
\right|}\right).
\]
Finally, we get the solid angle function $\Omega_i(\bx)$ of the facet $\mathcal{F}_\mathcal{P}^i$ for any $\bx\in\mathbf{R}^3:$
\[
\Omega_i(\bx)=2\pi-\sum_{j=1}^{n_i} \arccos\left(\frac{\left\langle \overrightarrow{XV}_{i_{j-1}}\times\overrightarrow{XV}_{i_j},
\overrightarrow{XV}_{i_{j}}\times\overrightarrow{XV}_{i_{j+1}} \right\rangle}{\left|\overrightarrow{XV}_{i_{j-1}}\times\overrightarrow{XV}_{i_j}\right|
\left|\overrightarrow{XV}_{i_{j}}\times\overrightarrow{XV}_{i_{j+1}}\right|}\right),
\]
\item We can summarize our results in the following
\begin{theorem}
Let us given a solid angle $\alpha$ ($0<\alpha<2\pi$) and a convex polyhedron $\mathcal{P}$ by its data structure
and its set of inequality. Then its isoptic surface can be determined by the equation
\begin{equation}
\alpha=\sum_{i=1}^m \mathbb{I}^{i}_{\mathcal{P}}(\bx)\Omega_i(\bx).
\label{eq1}
\notag
\end{equation}
\end{theorem}
\end{enumerate}
The results and computations will be demonstrated in the following subsection through the computation related to the regular tetrahedron and along with some figures.
\begin{remark}
\begin{enumerate}
	\item The algorithm can be easily extended for non-closed directed surfaces \textit{e. g.} for subdivision surfaces.
	\item If we have a \textbf{convex} polyhedron, than projecting its whole surface to the unit sphere,
	we obtain a double coverage (double solid angle) of the given polyhedron,
	therefore the algorithm can be changed \textit{i.e.} it is not necessary to determine the visible faces.
	In this case the isoptic surfaces are determined by the following \textbf{implicit} equation:
	\begin{equation}
	\alpha=\frac{1}{2}\sum_{i=1}^m \Omega_i(\bx).
	\label{eq2}
	\notag
\end{equation}
\end{enumerate}
\end{remark}
\subsection{Computations of the isoptic surface of a given regular tetrahedron}
Following the steps of the above described algorithm, we will calculate the compound isoptic surface of a given regular tetrahedron $\mathcal{T}$ whose data structure is determined by its vertices and faces where the faces
$\mathcal{F}_\mathcal{P}^i$ are given by their clockwise ordered vertices:
\begin{equation}
\begin{gathered}
 A_1=\left(0,0,\sqrt{\frac{2}{3}}-\frac{1}{2 \sqrt{6}}\right),  A_2=\left(-\frac{1}{2
   \sqrt{3}},-\frac{1}{2},-\frac{1}{2 \sqrt{6}}\right), \\ A_3=\left(-\frac{1}{2
   \sqrt{3}},\frac{1}{2},-\frac{1}{2 \sqrt{6}}\right),  A_4=
   \left(\frac{1}{\sqrt{3}},0,-\frac{1}{2 \sqrt{6}}\right),\\
   \left\{\mathcal{F}_\mathcal{T}^1,\left\{A_2,A_3,A_4\right\}\right\}, \left\{\mathcal{F}_\mathcal{T}^2,\left\{A_3,A_2,A_1\right\}\right\}, \\
 \left\{\mathcal{F}_\mathcal{T}^3,\left\{A_4,A_1,A_2\right\}\right\},  \left\{\mathcal{F}_\mathcal{T}^4,\left\{A_1,A_4,A_3\right\}\right\}.
\end{gathered} \notag
\end{equation}
The tetrahedron is also given by its system of inequality:
\[\left(
\begin{array}{ccc}
0&0&-4\sqrt{3}\\
-8\sqrt{6}&0&4\sqrt{3}\\
4\sqrt{6}&-12\sqrt{3}&4\sqrt{3}\\
4\sqrt{6}&12\sqrt{3}&4\sqrt{3}
\end{array}\right)\left(\begin{array}{c}
x\\
y\\
z
\end{array}\right)\leq
\left(\begin{array}{c}
\sqrt{2}\\
3\sqrt{2}\\
3\sqrt{2}\\
3\sqrt{2}
\end{array}\right).
\]
Since the final formula of the isoptic surface is too long to appear in print even for this simple example, we choose $i=1$ and derive only $\mathbb{I}^{1}_{\mathcal{P}}(\bx)$.  
\[
\mathbb{I}^{1}_{\mathcal{P}}(x,y,z)=\left\{
\begin{array}{lr}
	1 & -4\sqrt{3}z>\sqrt{2}       \\
	0 & -4\sqrt{3}z\leq\sqrt{2}
\end{array}
\right..
\]

For the resulted surface, see Fig. 3. (left), and in Fig. 3. (right) we describe a part of the tetrahedral isoptic surface of angle $\alpha=\pi/8$. 
The surface $\phi_{1,2}$ is derived by the solid angles belonging to the faces $\mathcal{F}_{\mathcal{T}}^3$ and $\mathcal{F}_{\mathcal{T}}^4$.
The isoptic surface $\phi_{1}$ is derived by the solid angles belonging to the faces $\mathcal{F}_{\mathcal{T}}^2$, $\mathcal{F}_{\mathcal{T}}^3$ 
and $\mathcal{F}_{\mathcal{T}}^4$ and
The isoptik surface $\phi_{2}$ is derived by the solid angles belonging to the faces $\mathcal{F}_{\mathcal{T}}^1$, $\mathcal{F}_{\mathcal{T}}^3$
and $\mathcal{F}_{\mathcal{T}}^4$.
\begin{figure}[H]
	\centering
		\includegraphics[width=6cm]{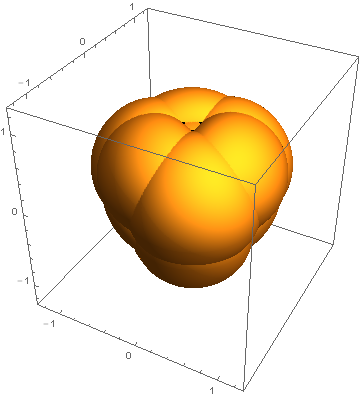}\includegraphics[width=6cm]{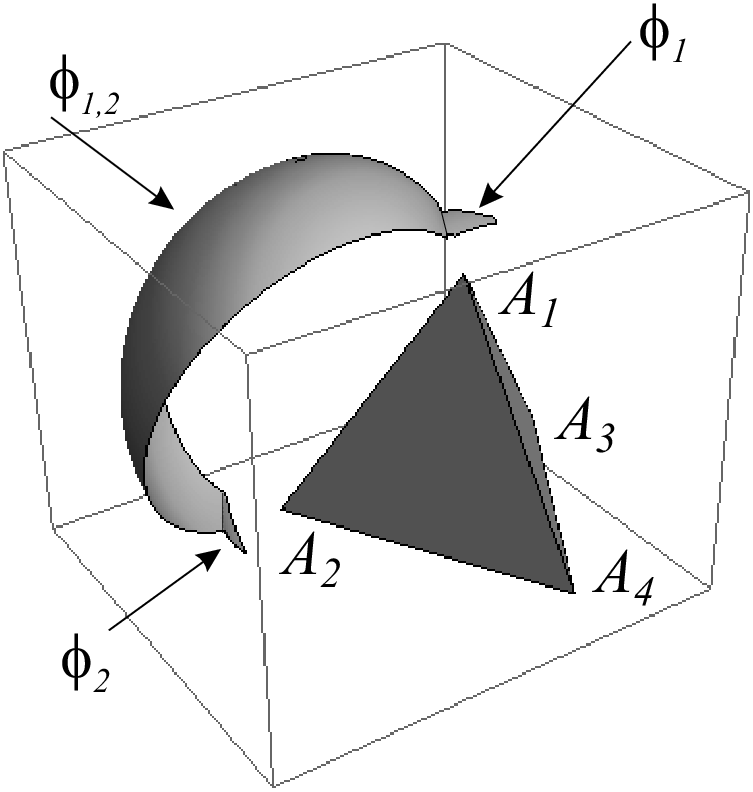}
	\caption{Isoptic surface of the regular tetrahedron for $\alpha=\pi/8$. A part of the tetrahedral isoptic surface e.g.
	The surface $\phi_{1,2}$ is derived by the solid angles belonging to the faces $\mathcal{F}_{\mathcal{T}}^3$ 
	and $\mathcal{F}_{\mathcal{T}}^4$.}
\end{figure}
\subsection{Isoptic surfaces to regular polyhedra and some Archimedean polyhedra}
\begin{figure}[H]
	\centering
		\includegraphics[width=6cm]{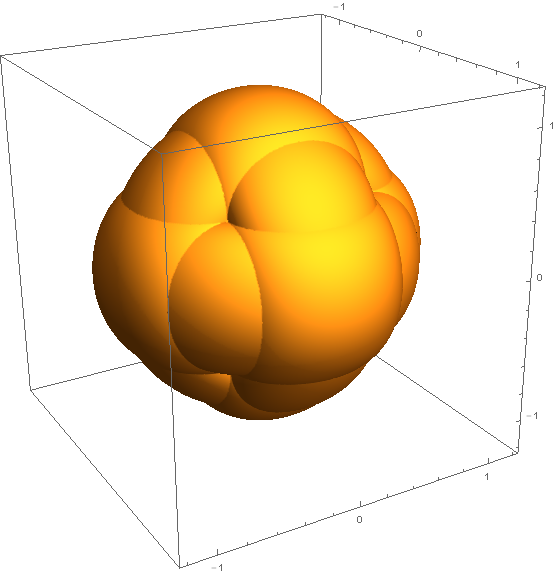}\includegraphics[width=6cm]{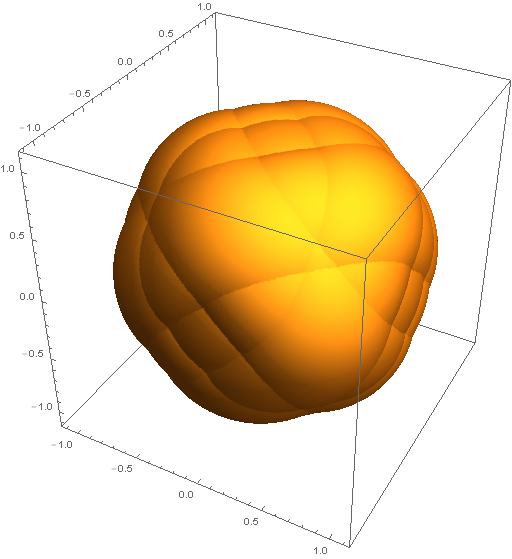}
	\caption{Isoptic surface of the cube $\alpha=\pi/2$ (left),\newline Isoptic surface of the regular octahedron for $\alpha=\pi/7$ (right)}
\end{figure}

\begin{figure}[H]
	\centering
		\includegraphics[width=6cm]{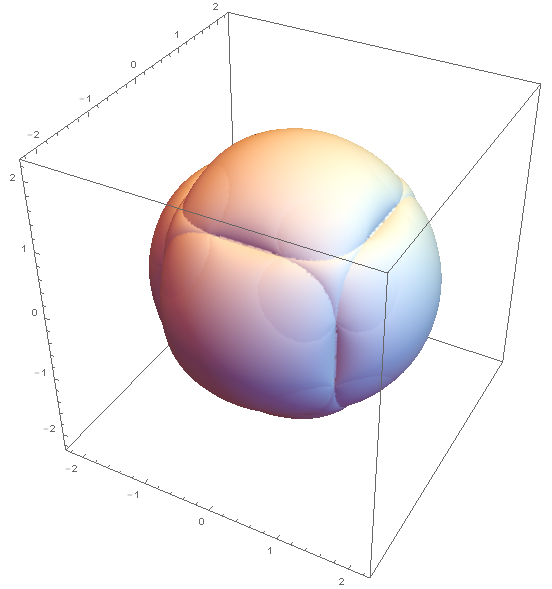} \includegraphics[width=6cm]{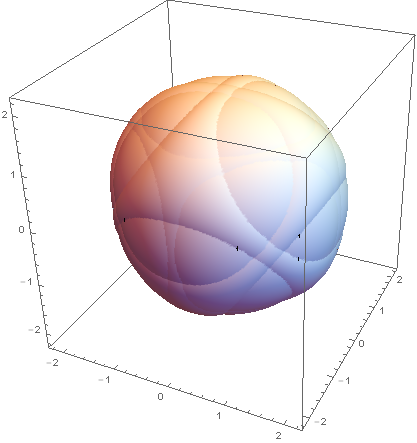}
	\caption{Isoptic surface of the truncated cube $\alpha=\pi$ (left)\newline Isoptic surface of the truncated octahedron $\alpha=2\pi/3$}
\end{figure}


\end{document}